\documentclass{amsproc}
\usepackage[margin=1.2in,nomarginpar]{geometry}
\usepackage{amssymb}
\usepackage{amsmath}
\usepackage{amsthm}

\usepackage{url}
\usepackage{xcolor}
\usepackage{amsmath,amssymb,amsthm,amsfonts,longtable}
\usepackage[english]{babel}
\usepackage{tikz-cd}
\usetikzlibrary{cd}
\usetikzlibrary{decorations.markings}
\tikzset{negated/.style={
		decoration={markings,
			mark= at position 0.5 with {
				\node[transform shape] (tempnode) {$\times$};
			}
		},
		postaction={decorate}
	}
}
\usepackage{array}
\usepackage[colorlinks,citecolor=red,urlcolor=blue,bookmarks=false,hypertexnames=true]{hyperref} 

\newtheorem{theorem}{Theorem}
\newtheorem{corollary}[theorem]{Corollary}
\newtheorem{proposition}[theorem]{Proposition}
\newtheorem{lemma}[theorem]{Lemma}
\newtheorem{remark}[theorem]{Remark}

\newcommand{\restr}{\mathord\downarrow} 
\newcommand{\ind}{\mathord\uparrow} 
\newtheorem{example}[theorem]{Example}
\newtheorem{definition}{Definition}[subsection]
\newcommand{\Irr}{\textnormal{Irr}}
\newcommand{\cd}{\textnormal{cd}}
\newcommand{\nl}{\textnormal{nl}}
\newcommand{\lin}{\textnormal{lin}}

\newcommand{\Core}{\textnormal{Core}}
\newcommand{\inertiagroup}{\textnormal{I}}
\newcommand{\cod}{\textnormal{cod}}
\newcommand{\CM}{\textnormal{CM}}
\newcommand{\Kern}{\textnormal{Kern}}

\title[Faithful quasi-permutation representations]{On the relation of character codegrees and the minimal faithful \\quasi-permutation representation degree of $p$-groups}
\author{Sunil Kumar Prajapati$^*$}
\address{Indian Institute of Technology, Bhubaneswar, Arugul Campus, Jatni, Khurda-752050, India.}
\email{skprajapati@iitbbs.ac.in}
\author{Ayush Udeep}
\address{Indian Institute of Technology, Bhubaneswar, Arugul Campus, Jatni, Khurda-752050, India.}
\email{udeepayush@gmail.com}
\thanks{$^{\textbf{*}}$ Corresponding author.}
\begin{document}
	\maketitle

	\begin{abstract}
		For a finite group $G$, we denote by $c(G)$, the minimal degree of faithful representation of $G$ by quasi-permutation matrices over the complex field $\mathbb{C}$. For an irreducible character $\chi$ of $G$, the codegree of $\chi$ is defined as $\cod(\chi) = |G/ \ker(\chi)|/ \chi(1)$. In this article, we establish equality between $c(G)$ and a $\mathbb{Q}_{\geq 0}$-sum of codegrees of some irreducible characters of a non-abelian $p$-group $G$ of odd order.
		We also study the relation between $c(G)$ and irreducible character codegrees for various classes of non-abelian $p$-groups, such as, $p$-groups with cyclic center, maximal class $p$-groups, GVZ $p$-groups, and  others.
	\end{abstract}
	
	\section{Introduction}
	Throughout this article, $G$ is a finite group and $p$ is a prime. The minimal faithful permutation degree $\mu(G)$ of $G$ is the least positive integer $n$ such that $G$ is isomorphic to some subgroup of $S_{n}$. In a parallel direction to the definition of $\mu(G)$, another degree $c(G)$ is defined as the minimal degree of a faithful representation of $G$ by complex quasi-permutation matrices (square matrices over complex field with non-negative integral trace) (see \cite{BGHS}). Since every permutation matrix is a quasi-permutation matrix, it is easy to see that $c(G)\leq \mu(G)$.  In fact, Behravesh and Ghaffarzadeh \cite[Theorem 3.2]{BG} proved that if $G$ is a finite $p$-group of odd order, then $c(G)=\mu(G)$.
	Several researchers have studied $\mu(G)$ and $c(G)$ extensively in the past (see \cite{HB, BG, BGHS, GA, DLJ, DW}). Behravesh \cite[Theorem 3.6]{HB} gave an algorithm for the computation of $c(G)$. 
	In \cite{BG}, Behravesh and Ghaffarzadeh improved the said algorithm and proved the following.
\begin{lemma} \textnormal{\cite[Lemma 2.2]{BG}} \label{lemma:c(G)Algorithm}
	Let $G$ be a finite group.  Let $X \subset 
	\Irr(G)$ be such that $\cap_{\chi \in X} \ker (\chi)= 1$ 
	and $\cap_{\chi \in Y} \ker (\chi) \neq 1$ 
	for every proper subset $Y$ of $X$. 
	Let $\xi_X = \sum_{\chi \in X} \left[ \sum_{\sigma \in \Gamma(\chi)}
	\chi^{\sigma}  \right]$ and let $m(\xi_X)$ be the absolute value of 
	the minimum value that $\xi_X$ takes over $G$.
	Then $$c(G) = \min \{\xi_X(1) + m(\xi_X) \; | \; X \subset \Irr(G) 
	\text{\ satisfying\ the\ above\ property} \}.$$
\end{lemma}

\noindent 
We identify $X_G\subset \Irr(G)$ with a minimal 
faithful quasi-permutation representation of $G$ if 
\begin{equation}\label{eq:X_G} 
	\bigcap_{\chi \in X_G} \ker (\chi) = 1 \text{ and } 	
	\bigcap_{\chi \in Y} \ker (\chi) \neq 1 
	\text{ for every } Y \subset X_G
\end{equation}
and $c(G) = \xi_{X_G}(1) + m(\xi_{X_G})$. Let $\mathcal{S}$ be the collection of all minimal faithful quasi-permutation representations of $G$ satisfying \eqref{eq:X_G}. \\

	Now, let $\chi$ be an irreducible character of $G$ and $\Gamma(\chi)$ be the Galois group of $\mathbb{Q}(\chi)$ over $\mathbb{Q}$. The codegree of $\chi$ is defined as $\cod(\chi) = |G/ \ker(\chi)|/ \chi(1)$, and the set of codegrees of all irreducible characters of $G$ is denoted by $\cod(G)$ (see \cite{CLcodegree, CLmaximal, DLcodegree,QWWcodegree}). When $\chi$ is linear, $\cod(\chi) = |G/ \ker(\chi)|$, and it is easy to observe that $$\sum_{\sigma \in \Gamma(\chi)} \chi^{\sigma}(1) = \chi(1)[\mathbb{Q}(\chi): \mathbb{Q}] = \chi(1)|\Gamma(\chi)| = |\Gamma(\chi)| = \phi(\cod(\chi)).$$ Therefore by Lemma \ref{lemma:c(G)Algorithm}, $c(G)$ is  related to the codegrees of $\chi\in X_G$. 

\begin{example} \label{example:intro1}
	\textnormal{Let
	$G$ be an extraspecial $p$-group of order $p^3$ ($p \geq 3$). Then $Z(G) \cong C_p$ and $\cd(G) = \{ 1, p \}$. Observe that there exists an elementary abelian subgroup of index $p$ in $G$. Then from \cite[Proposition 19]{SAcyclic}, $c(G) = p^2$. Let $X_G \in \mathcal{S}$. It is easy to see that $X_G = \{\chi \}$, where $\chi$ is any faithful irreducible character of $G$. Then $\cod(\chi) = |G/ \ker(\chi)|/ \chi(1) = p^2$, which implies that $c(G) = \sum_{\chi \in X_G}\cod(\chi)$.}
\end{example}
\begin{example} \label{example:intro2}
	\textnormal{Suppose
		\[ G =  \langle x, y, z: x^{p^2} = y^p = z^p = 1, xy = yx^{p+1}, xz = zxy, yz = zy \rangle, \]
		which is a $p$-group of order $p^4$ ($p \geq 3$) with $Z(G) = \langle x^p \rangle \cong C_p$ and $\cd(G) = \{ 1, p \}$. Observe that $\langle x^p, y, z \rangle$ is an elementary abelian subgroup of index $p$ in $G$. Then from \cite[Proposition 19]{SAcyclic}, $c(G) = p^2$. Let $X_G \in \mathcal{S}$. Then $X_G = \{ \chi \}$, where $\chi$ is any faithful irreducible character of $G$. 
		 Hence, $\cod(\chi) = |G/ \ker(\chi)|/ \chi(1) = p^3$, which implies that $c(G) < \sum_{\chi \in X_G}\cod(\chi)$.
	}
\end{example}
\noindent	In this article, we study the relation between $c(G)$ and character codegrees of a non-abelian $p$-group $G$. To the best of our knowledge, this relation has not been studied in the past.\\
	
	Suppose $G$ is a non-abelian $p$-group $(p\geq 3)$ and let $X_G \in \mathcal{S}$. From Lemma \ref{thm:ford}, for all  $\chi \in X_G$, there exists a subgroup $H_{\chi}$ of $G$ and $\lambda_{\chi} \in \lin(H_{\chi})$ such that $\chi = \lambda_{\chi}\ind_{H_{\chi}}^{G}$ and $\mathbb{Q}(\chi) = \mathbb{Q}(\lambda)$.
	Define $a_{\chi} = |\ker(\lambda_{\chi}):\ker(\chi)|$.
	Then in Theorem \ref{lemma:generalcodegreerelation}, we prove that 
		\begin{equation*} \label{eq:codrelationNotbestboundIntro}
				c(G) =  \sum_{\chi \in X_{G} \cap \lin(G)}\cod(\chi) + \sum_{\chi \in X_{G}\cap \nl(G)} \frac{\chi(1)}{a_{\chi}} \cod(\chi), \text{ where } a_{\chi} \text{ divides } \cod(\chi).
			\end{equation*}
	 Since $a_{\chi} \geq 1$, for each $\chi\in X_G\cap \nl(G)$, we get
	 \begin{equation} \label{eq:codgeneral}
	 	c(G) \leq \sum_{\chi \in X_{G}} \chi(1)\cod(\chi).
	 \end{equation}
	  Now, suppose $\ker(\lambda_{\chi}) \text{ is not normal in } G$, for each $\chi \in X_G \cap \nl(G)$. Then $a_{\chi} \geq p$, for each $\chi \in X_G \cap \nl(G)$, and we get
	  \begin{equation} \label{eq:introbestbound}
	  	c(G) \leq \sum_{\eta \in X_{G}\cap \lin(G)} \cod(\eta) + \sum_{\chi \in X_{G}\cap \nl(G)} \frac{1}{p}\chi(1) \cod(\chi),
	  \end{equation}
	  which is an improved bound from Inequality \eqref{eq:codgeneral}. Moreover, there exists a $p$-group ($p\geq 3$) with cyclic center for which \eqref{eq:introbestbound} is an equality (see Proposition \ref{prop:codnMcyclic}).
	   In many instances, we have $a_{\chi} = p$ for all $\chi \in X_G \cap \nl(G)$ (for example, in a non-abelian core-$p$ $p$-group ($p\geq 3$) with cyclic center).
	This motivates us to define a class $\mathcal{T}$ of non-abelian $p$-groups ($p\geq 3$) in the following way.
	\begin{align*}
	\text{\bf Hypothesis I: }	\text{ For every } X_{G} \in \mathcal{S}, \ker(\lambda_{\chi}) \text{ is NOT normal in $G$}, \text{ for all }  \chi \in X_{G} \cap \nl(G).
	\end{align*}
	 We define $\mathcal{T}$ as follows:
	\[ \mathcal{T} = \{ G : |G| = p^n~ (p\geq 3) \text{ and } G \text{ satisfies Hypothesis I} \}. \]

	\noindent We obtain many classes of non-abelian $p$-groups belonging in $\mathcal{T}$, namely, $p$-groups with a cyclic center (see Lemma \ref{lemma:cycliccenterinT}), $p$-groups $(p\geq 3)$ with order $\leq p^4$ (see Theorem \ref{thm:codlessthanp4}), some class of core-$p$ $p$-groups (see Subsection \ref{subsec:corep}), $p$-groups where kernels of nonlinear irreducible characters form a chain (see Subsection \ref{subsec:chainofkernels}), and GVZ $p$-groups (see Subsection \ref{subsec:GVZgroups}).\\
	
	 We also establish equality between $c(G)$ and  the sum of codegrees of $\chi \in X_G$ for various classes of $p$-groups in $\mathcal{T}$.
	 A $p$-group $G$ is called a core-$p$ $p$-group if $|H/ \Core_{G}(H)| \leq p$ for every subgroup $H$ of $G$, where $\Core_{G}(H)$ is the core of $H$ in $G$ (see \cite{corep1995}). In Theorem \ref{thm:kerlambdanotnormal}, we prove that if $G$ is a core-$p$ $p$-group $(p\geq 3)$ belonging to $\mathcal{T}$ then $c(G) = \sum_{\chi \in X_{G}}\cod(\chi)$. Further, if $Z(G)$ is cyclic, then $c(G) = |G|/p$.
	 
	 A group $G$ is called a GVZ-group if $\chi(g) = 0$, for all $g\in G\setminus Z(\chi)$ and all $\chi \in \nl(G)$.
	In Theorem \ref{thm:GVZcod}, we prove the following.
		\begin{theorem} \label{thm:GVZcod}
		Let $G$ be a GVZ $p$-group.
		\begin{enumerate}
			\item [\rmfamily (i)] Let $p\geq 3$. For each $\chi \in \Irr(G)$, choose $\lambda_{\chi} \in \lin(H_{\chi})$, for some $H_{\chi} \leq G$ such that $\chi = \lambda_{\chi}\ind_{H_{\chi}}^{G}$ and $\mathbb{Q}(\chi) = \mathbb{Q}(\lambda)$. Then $|\ker(\lambda_{\chi})| = \chi(1)|\ker(\chi)|$. Further, $ G\in \mathcal{T}$.  
			\item [\rmfamily (ii)] If $p$ is any prime and $X_G \in \mathcal{S}$, then $c(G) = \sum_{\chi \in X_{G}} \cod(\chi)$.
		\end{enumerate}
	\end{theorem}
\noindent With the help of Theorem \ref{thm:GVZcod}, we prove that if $G$ is a GVZ $p$-group with cyclic center, then $c(G) = |G/Z(G)|^{1/2}|Z(G)|$ (see Corollary \ref{cor:GVZcyclic}).

We also study a special sub-class of GVZ $p$-groups called $\CM_{p-1}$ $p$-groups. 
A group $G$ is called a $\CM_{n}$-group if every normal subgroup of $G$ appears as the kernel of at most $n$ irreducible characters of $G$ (see \cite{BL}). In Theorem \ref{thm:CMgroup}, we prove the following result for $\CM_{p-1}$-groups.

	\begin{theorem} \label{thm:CMgroup}
	Let $G$ be a $\CM_{p-1}$ $p$-group. Then
	\begin{enumerate}
		\item [(i)] $\cod(G) = \{ 1, p\chi(1)~ |~ 1_{G} \neq \chi \in \Irr(G) \}$.
		\item [(ii)] Suppose $X_G \in \mathcal{S}$. Then
		\[ c(G) = p \sum_{\chi \in X_{G}} \chi(1) = \sum_{\chi \in X_{G}} \cod(\chi). \]
		\item [(iii)] Suppose $\cd(G) = \{ 1 = d_{0}, d_{1}, \ldots, d_{t} \}$. If $X_G, Y_G \in \mathcal{S}$ such that $Y_{G} \neq X_{G}$, then $|\Irr_{d_{i}}(G) \cap Y_{G}| = |\Irr_{d_{i}}(G) \cap X_{G}|$, for each $i~ ( 0\leq i \leq t)$.
	\end{enumerate}
\end{theorem}

	 	We begin with Section \ref{sec:notations}, where we mention the notations and preliminary results required in this article.
	In Subsection \ref{subsec:corep}, we study classes of core-$p$ $p$-groups that belong to $\mathcal{T}$.
 In Subsection \ref{subsec:chainofkernels}, we study the $p$-groups where kernels of nonlinear irreducible characters form a chain. 
In Subsection \ref{subsec:GVZgroups}, we study GVZ $p$-groups.
	   We conclude our article by discussing results and posing a few questions to our readers in Section \ref{sec:results}.

	\section{Notations and Preliminaries} \label{sec:notations}
	
	For a finite group $G$, we use the following notations throughout the article.\\
	\begin{tabular}{cl}
		$d(G)$ & the minimal number of generators of $G$\\
		$G'$ & the commutator subgroup of $G$\\
		$\exp(G)$ & the exponent of $G$\\
		$\Core_{G}(H)$ & the core of $H$ in $G$, for $H\leq G$\\
		$\Irr(G)$ & the set of irreducible complex characters of $G$\\
		$\Irr_{t}(G)$ & the set of irreducible complex characters of $G$ of degree $t$\\
		$\lin(G)$ & the set of linear characters of $G$\\
		$\nl(G)$ & the set of non-linear irreducible characters of $G$\\
		$\cd(G)$ & the character degree set of $G$, i.e., $\{ \chi(1) ~|~ \chi \in \Irr(G) \}$\\
		$Z(\chi)$ & $\{ g\in G~|~ |\chi(g)| = \chi(1) \}$, for $\chi \in \Irr(G)$\\
		$b(G)$ & $\max \{ \chi(1) ~|~ \chi \in \Irr(G) \}$\\
		$\mathbb{Q}(\chi)$ &  the field obtained by adjoining the values $\chi(g)$ to $\mathbb{Q}$, for $\chi \in \Irr(G)$  and all
		$g \in G$\\
		$\Gamma(\chi)$ & the Galois group of $\mathbb{Q}(\chi)$ over $\mathbb{Q}$, for $\chi \in \Irr(G)$\\
%
		$\phi(n)$ & the Euler phi function\\
		$\omega_{n}$ & a primitive $n^{th}$ root of unity\\
		$\mathcal{S}$ & the collection of all minimal faithful quasi-permutation representations of $G$ satisfying \eqref{eq:X_G}\\
		$\mathbb{Q}_{\geq 0}$ & the set of non-negative rational numbers
	\end{tabular}
	\\
	
	Now we summarize some results that 
	we use throughout the article. Throughout this article, unless otherwise 
	indicated, $p$ is an arbitrary prime.\\
	
	For irreducible characters of finite $p$-groups, Ford has proved the following result in \cite{FORD}.
	\begin{lemma}\textnormal{\cite[Theorem 1]{FORD}} \label{thm:ford}
		Let $G$ be a $p$-group $(p\geq 3)$ and let $\chi$ be an irreducible complex character 
		of $G$. Then there exists a linear character $\lambda$ on a subgroup 
		$H$ of $G$ such that 
		$\lambda\ind_{H}^{G} = \chi$ and $\mathbb{Q}(\lambda) = \mathbb{Q}(\chi)$.
	\end{lemma}
	
	\begin{remark} \label{remark:ford}
		\textnormal{ Let $G$ be a $p$-group $(p\geq 3)$ and let $\chi\in \Irr(G)$. Throughout the article, we use Lemma \ref{thm:ford} extensively. Therefore, to avoid repetition, let $H_{\chi}$ denote a subgroup of $G$ and let $\lambda_{\chi}$ denote a linear character of $H_{\chi}$ such that $\lambda_{\chi}\ind_{H_{\chi}}^{G} = \chi$ and $\mathbb{Q}(\lambda_{\chi}) = \mathbb{Q}(\chi)$, unless stated otherwise. }
	\end{remark}
	
	\noindent Let $\chi, ~\psi \in \Irr(G)$. 
	We say that $\chi$ and $\psi$ are Galois conjugate over $\mathbb{Q}$ if there exists $\sigma \in \Gamma(\chi)$ such that $\chi^{\sigma}= \psi$. One can check that Galois conjugacy defines an equivalence relation on $\Irr(G)$. Moreover, if $\mathcal{C}$ denotes the equivalence class of $\chi$ with respect to Galois conjugacy over $\mathbb{Q}$, then $|\mathcal{C}|=| \mathbb{Q}(\chi) : \mathbb{Q} |$ (see \cite[Lemma 9.17]{I}).

		\begin{definition} \label{D1}
		Let $G$ be a finite group. 
		\begin{enumerate}
			\item [\rmfamily(i)] For $\psi\in \Irr(G)$, define $d(\psi)= |\Gamma(\psi)|\psi(1)$.
			\item [\rmfamily(ii)] For any complex character $\chi$ of $G$, define 
			\[  m(\chi)=
			\begin{cases}
				0 &\quad \text{ if } \chi(g) \geq 0 \text{ for all } g\in G,\\
				-\min \left\{ \sum\limits_{\sigma\in \Gamma(\chi)} \chi^{\sigma}(g): g\in G \right\} &\quad \text{ otherwise.} 
			\end{cases} \]
		\end{enumerate}
	\end{definition}

	\begin{lemma}\textnormal{\cite[Theorem 2.3]{GA}}\label{lemma:m(chi)}
		Let $G$ be a non-abelian $p$-group and let $X_{G}\in \mathcal{S}$. Then
		\begin{enumerate}
			\item [\rmfamily(i)] $|X_G| = d(Z(G))$.
			\item [\rmfamily(ii)] Let $\xi = \sum_{\chi \in X_G}  \sum_{\sigma \in \Gamma(\chi)}
			\chi^{\sigma} $. Then
			\[ m(\xi) = \frac{1}{p-1} \sum_{\chi \in X_{G}}
		 \sum_{\sigma \in \Gamma(\chi)} 
			\chi^{\sigma}(1). \]
		\end{enumerate}
	\end{lemma}
	
	\begin{remark} \label{remark:m(chi)}
		\textnormal{Let $G$ be a non-abelian $p$-group and suppose $X_G \in \mathcal{S}$. Let $\xi = \sum_{\chi \in X_G} \sum_{\sigma \in \Gamma(\chi)}
			\chi^{\sigma} $.
			\begin{enumerate}
				\item [(i)] From Lemma \ref{lemma:m(chi)}, $m(\xi) = \frac{1}{p-1}\xi(1)$. Then
				\[ c(G) = \xi(1) + m(\xi) = \xi(1) + \frac{1}{p-1}\xi(1) = \frac{p}{p-1}\sum_{\chi\in X_G}\chi(1)|\Gamma(\chi)| = \frac{p}{p-1}\sum_{\chi\in X_G} d(\chi). \]
					\item [(ii)]  Let $\chi \in X_G$. There exists a central element $z$ of order $p$ in $G$ (see the proof of \cite[Theorem 2.3]{GA}) such that 
				\[ \sum_{\sigma \in \Gamma(\chi)}
				\chi^{\sigma}(z) = -\frac{1}{p-1}\chi(1)|\Gamma(\chi)| \Rightarrow m(\chi) \geq \frac{1}{p-1}d(\chi). \]
				Hence, $\sum_{\chi \in X_G} m(\chi) \geq \sum_{\chi \in X_G} \frac{1}{p-1}d(\chi) = m(\xi)$.
			\end{enumerate} }
	\end{remark}

\begin{lemma} \label{thm:mu(G)=c(G)}
	\textnormal{\cite[Theorem 3.2]{BG}}
	\label{thm:pgroup}
	If $G$ is a finite $p$-group $(p\geq 3)$, then $c(G) = \mu(G)$.
\end{lemma}

	\begin{lemma} \textnormal{\cite[{Theorem 7}]{SAcyclic}} 
		\label{lemma:normallymonomial}
		Let $G$ be a normally monomial $p$-group with cyclic center.
		Suppose $A$ is an  abelian normal subgroup of maximum order in $G$. Then 
		\[ (\max \cd(G))|Z(G)| \text{ divides } c(G) \text{ and } c(G) \text{ divides } (\max \cd(G))\exp(A). \]
	\end{lemma}
	
	
		\section{Results} \label{sec:discussions}
We start with the following lemma.
	
		\begin{lemma} \label{lemma:lambdaisnormal}
		Let $G$ be a $p$-group $(p\geq 3)$ and $\chi \in \nl(G)$. Suppose, for each $i=1,2$, $\chi = \lambda_i\ind_{H_i}^{G}$, for some $H_i \leq G$ and for some $\lambda_i \in \lin(H_i)$ such that $\mathbb{Q}(\chi) = \mathbb{Q}(\lambda_i)$. Then $\ker(\lambda_1)$ is normal in $G$ if and only if $\ker(\lambda_{2})$ is normal in $G$.
	\end{lemma}
	\noindent \emph{Proof.} Let $\ker(\lambda_1)$ be a normal subgroup of $G$. Then $\ker(\lambda_1) = \ker(\chi)$. Now, let $\gamma = \sum_{\sigma\in \Gamma(\chi)} \chi^{\sigma}$. Then 
	$\gamma(1) = \chi(1)|\Gamma(\chi)| = \chi(1) |\Gamma(\lambda_1)| = \chi(1) \phi(|H_1/\ker(\lambda_1)|)$. Similarly, $\gamma(1) = \chi(1) |\Gamma(\lambda_2)| = \chi(1) \phi(|H_2/\ker(\lambda_2)|)$. Since $|H_1| = |H_{2}|$,  we get $|\ker(\chi)| = |\ker(\lambda_1)| = |\ker(\lambda_{2})|$. Therefore, $\ker(\lambda_{2})$ is also normal in $G$.\\ In the same way, if $\ker(\lambda_1)$ is normal in $G$, then $\ker(\lambda)$ is normal in $G$. \qed

	\begin{lemma} \label{lemma:cycliccenterinT}
		Let $G$ be a non-abelian $p$-group $(p\geq 3)$ with cyclic center. Then $G \in \mathcal{T}$.
	\end{lemma}
	\noindent \emph{Proof.} From Lemma \ref{lemma:m(chi)}, we get $|X_{G}| = 1$, where $X_G \in \mathcal{S}$. Moreover, $X_{G} \cap \lin(G) = \emptyset$. Now, suppose $X_{G} = \{ \chi \}$, for some faithful irreducible character $\chi$ of $G$. Let $H_{\chi} \leq G$ and $\lambda_{\chi} \in \lin(H_{\chi})$ be as defined in Remark \ref{remark:ford}.\\
	{\bf Claim:} $\ker(\lambda_{\chi}) \neq 1$. On the contrary, suppose that $\ker(\lambda_{\chi}) = 1$. Then $d(\chi) = \chi(1)|\Gamma(\chi)| = \chi(1) |\Gamma(\lambda_{\chi})| = \chi(1)\phi(|H_{\chi}/\ker(\lambda_{\chi})|) = |G/H_{\chi}|\phi(|H_{\chi}|) = \phi(|G|)$. Then from Remark \ref{remark:m(chi)}, we get $c(G) = |G|$, which is a contradiction since $G$ is a non-abelian $p$-group ($p\geq 3$). Hence, $\ker(\lambda_{\chi}) \neq 1$.\\
	 Then $1 = |\ker(\chi)| < |\ker(\lambda_{\chi})|$. Therefore, $\ker(\lambda_{\chi})$ is not normal in $G$, and hence, $G \in \mathcal{T}$. \qed\\
	
	In Lemma \ref{lemma:codlinearchar}, we prove that the codegree of a nontrivial linear character of a $p$-group is a multiple of the order of its Galois conjugacy class.
	
		\begin{lemma} \label{lemma:codlinearchar}
		Let $G$ be a $p$-group and $1_G \neq \eta \in \lin(G)$. Then 
		\[ \cod(\eta) = \frac{p}{p-1} d(\eta). \]
	\end{lemma}
	\noindent \emph{Proof.} Let $G$ be a $p$-group and $1_G \neq \eta \in \lin(G)$.  Suppose $\cod(\eta) = |G/\ker(\eta)| = p^{b}$ (say). Then $d(\eta) = \eta(1) |\Gamma(\eta)| = \phi\left( |G/\ker(\eta)| \right) = \phi(p^{b})$. From \cite[Lemma 4.5]{HB}, $m(\eta) = p^{b-1} = \frac{1}{p-1}d(\eta)$. Then \[ d(\eta) + m(\eta) = \frac{p}{p-1}d(\eta) = p^b = \cod(\eta). \eqno \qed \]

	\begin{theorem}\label{thm:codlessthanp4}
		Let $G$ be a non-abelian $p$-group of order $\leq p^4 ~ (p\geq 3)$. Then $G \in \mathcal{T}$. Further, let $X_G \in \mathcal{S}$. Then
		\begin{enumerate}
			\item [(i)] if $|G| = p^3$, then $c(G) = \sum_{\chi \in X_G} \cod(\chi)$, and
			\item [(iii)] if $|G| = p^4$, then $c(G) \leq \sum_{\chi \in X_G} \cod(\chi)$. Moreover, if $|Z(G)| = p^2$, then $c(G) = \sum_{\chi \in X_G} \cod(\chi)$.
		\end{enumerate}
	\end{theorem}
	
\noindent \emph{Proof.} If $G$ is a group of order $p^3$ $(p\geq 3)$, then $Z(G) \cong C_p$ and hence from Lemma \ref{lemma:cycliccenterinT}, $G \in \mathcal{T}$. Further, from Example \ref{example:intro1}, we get $c(G) = \sum_{\chi \in X_G} \cod(\chi)$, for any $X_G \in \mathcal{S}$. \\
Now let $|G| = p^4$ ($p\geq 3$) and suppose $X_G \in \mathcal{S}$. Then $|Z(G)| \leq p^2$ and $\cd(G) = \{ 1, p \}$. \\
\noindent {\bf Case I ($|Z(G)| = p$):} From Lemma \ref{lemma:cycliccenterinT}, we have $G \in \mathcal{T}$. From \cite[Section 4.4]{RJ}, $\exp(G) \leq p^2$. Then from \cite[Corollary 6]{SAcyclic}, we get $c(G) = p^2$ or $p^3$. From Lemma \ref{lemma:m(chi)}, $|X_G| = 1$. Let $X_{G} = \{ \chi \}$, for some faithful irreducible character $\chi$ of $G$.  Then $\cod(\chi) = p^3$ and thus, $c(G) \leq \cod(\chi)$.\\
	{\bf Case II ($|Z(G)| = p^2$):} If $Z(G)$ is cyclic, then from Lemma \ref{lemma:cycliccenterinT}, $G \in \mathcal{T}$ and from \cite[Corollary 6]{SAcyclic}, $c(G) = p^3$. Now let $X_{G} = \{ \chi \} \in \mathcal{S}$. Then $\ker(\chi) = 1$ and $\cod(\chi) = p^3$. Thus, we get $c(G) = \cod(\chi)$.\\
	 Now suppose $Z(G) \cong C_p \times C_p$. Let $X_G \in \mathcal{S}$ and suppose $\chi \in X_G \cap \nl(G)$.
	   Then $\chi(1) = p = |G/Z(G)|^{1/2}$, and thus from \cite[Corollary 2.30]{I}, $\chi(g) = 0$ for all $g\in G\setminus Z(G)$. Here $\ker(\chi) \subset Z(G)$, and so $|\ker(\chi)| = p$.
	   Let $H_{\chi} \leq G$ and $\lambda_{\chi} \in \lin(H_{\chi})$ be as defined in Remark \ref{remark:ford}. Since $Z(G) \cong C_p \times C_p$, we get $\mathbb{Q}(\chi) = \mathbb{Q}(\omega_{p})$. Then $\mathbb{Q}(\lambda_{\chi}) = \mathbb{Q}(\chi) = \mathbb{Q}(\omega_{p})$ and so, $|H_{\chi}/\ker(\lambda_{\chi})| = p$. This implies that $|\ker(\lambda_{\chi})| = p^2$.
	     Since $\ker(\chi) = \Core_{G}(\ker(\lambda_{\chi}))$, we deduce that $\ker(\lambda_{\chi})$ is not normal in $G$. Therefore, $G \in \mathcal{T}$. \\
	 Now, for all $\chi \in X_G \cap \nl(G)$, we get $\cod(\chi) = p^2$, and $\frac{p}{p-1}d(\chi) = \frac{p}{p-1}\chi(1)|\Gamma(\chi)| = \frac{p^2}{p-1}\phi(p)= p^2$.
 Then from Remark \ref{remark:m(chi)} and Lemma \ref{lemma:codlinearchar}, we get
	\[ c(G) = \frac{p}{p-1} \sum_{\chi \in X_{G}}d(\chi) = \sum_{\chi \in X_{G}} \cod(\chi). \eqno \qed \]

%
	
	Suppose $G$ is a non-abelian $p$-group $G$ $(p\geq 3)$ and let $X_G \in \mathcal{S}$. In Theorem \ref{lemma:generalcodegreerelation}, we establish equality between $c(G)$ and a $\mathbb{Q}_{\geq 0}$-sum of codegree of $\chi \in X_G$.
	\begin{theorem}
		\label{lemma:generalcodegreerelation}
		Let $G$ be a non-abelian $p$-group $(p\geq 3)$ and suppose $X_G \in \mathcal{S}$.
		\begin{enumerate}
			\item [(i)] For each $\chi \in X_G$, choose $\lambda_{\chi} \in \lin(H_{\chi})$, for some $H_{\chi} \leq G$ such that $\chi = \lambda_{\chi}\ind_{H_{\chi}}^{G}$ and $\mathbb{Q}(\chi) = \mathbb{Q}(\lambda)$. Then
			\begin{equation*} 
				c(G) =  \sum_{\chi \in X_{G} \cap \lin(G)}\cod(\chi) + \sum_{\chi \in X_{G}\cap \nl(G)} \frac{\chi(1)}{a_{\chi}} \cod(\chi), \text{ where } a_{\chi} = |\ker(\lambda_{\chi}) : \ker(\chi)|.
			\end{equation*}
			\item [(ii)] Further, if $G\in \mathcal{T}$, then
			 \begin{equation*} \label{eq:remarkcod}
				c(G) \leq \sum_{\eta \in X_{G}\cap \lin(G)} \cod(\eta) + \sum_{\chi \in X_{G}\cap \nl(G)} \frac{1}{p}\chi(1)\cod(\chi).
			\end{equation*}
		\end{enumerate} 
	\end{theorem}
	\noindent \emph{Proof.} Let $\chi \in X_G\cap \nl(G)$. Then
			\begin{align*}
				\left| \frac{H_{\chi}}{\ker(\lambda_{\chi})} \right| = \frac{|H_{\chi}|}{a_{\chi}|\ker(\chi)|} = \frac{|G|}{a_{\chi}|\ker(\chi)|\chi(1)} = \frac{\cod(\chi)}{a_{\chi}}, \text{ where } a_{\chi} = |\ker(\lambda_{\chi}) : \ker(\chi)|.
			\end{align*}
	So $d(\chi) = \chi(1)|\Gamma(\chi)| = \chi(1)|\Gamma(\lambda_{\chi})| = \chi(1)\phi\left( \left| \frac{H_{\chi}}{\ker(\lambda_{\chi})} \right|  \right) = \chi(1) \phi(\cod(\chi)/a_{\chi})$. Then 
	\begin{equation} \label{eq:gencodrelation}
		\frac{p}{p-1}d(\chi) = \frac{\chi(1)}{a_{\chi}}\cod(\chi), \text{ for each } \chi \in X_G \cap \nl(G).
	\end{equation} 
	Therefore, from Remark \ref{remark:m(chi)} and Lemma \ref{lemma:codlinearchar}, we get
	\[ 	c(G) = \frac{p}{p-1}\sum_{\chi \in X_{G}}  d(\chi) = \sum_{\chi \in X_{G} \cap \lin(G)}\cod(\chi) + \sum_{\chi \in X_{G}\cap \nl(G)} \frac{\chi(1)}{a_{\chi}} \cod(\chi).  \]
 Now, suppose $G \in \mathcal{T}$. Then $a_{\chi} \geq p$ for each $\chi \in X_G \cap \nl(G)$. 
	Hence from Remark \ref{remark:m(chi)} and Lemma \ref{lemma:codlinearchar}, we get
\[ 	c(G) \leq \sum_{\chi \in X_{G} \cap \lin(G)}\cod(\chi) + \sum_{\chi \in X_{G}\cap \nl(G)} \frac{\chi(1)}{p} \cod(\chi). \eqno \qed \]

	\begin{remark} \label{remark:codupperbound}
	\textnormal{We do {\bf NOT} have any example of a non-abelian $p$-group in which $a_{\chi} = 1$, for any $\chi \in X_G \cap \nl(G)$ and for an arbitrary $X_G \in \mathcal{S}$. This motivates us to  investigate the class $\mathcal{T}$, in which $a_{\chi} \geq p$ for all $\chi \in X_G \cap \nl(G)$. }
\end{remark}

\begin{proposition} \label{prop:newclassinT}
	Let $G$ be a non-abelian $p$-group $(p\geq 3)$ such that $\cd(G) = \{ 1, p^b \} ~(b\geq 1)$ and $c(G) = d(Z(G)) p^{b+1}$. Then $G \in \mathcal{T}$.
\end{proposition}
\noindent \emph{Proof.} Let $X_G \in \mathcal{S}$. For $\chi \in X_G$, let $H_{\chi} \leq G$ and $\lambda_{\chi} \in \lin(H_{\chi})$ be as defined in Remark \ref{remark:ford}. Let $|H_{\chi} / \ker(\lambda_{\chi})| = p^{t_{\chi}}$, for some $t_{\chi} \geq 1$.\\
{\bf Claim:} For $\chi \in X_G \cap \nl(G)$, we get $t_{\chi} = 1$.\\
 For $\chi \in X_G$, we have $d(\chi) = \chi(1)|\Gamma(\chi)| = \chi(1)|\Gamma(\lambda_{\chi})| = \chi(1)\phi(p^{t_{\chi}})$. Hence, from Remark \ref{remark:m(chi)}, we get
\[ c(G) = \frac{p}{p-1}\sum_{\chi \in X_{G}} d(\chi) =   \sum_{\chi \in X_{G}}\frac{p}{p-1}\chi(1) \phi(p^{t_{\chi}}) =  \sum_{\chi \in X_{G}} \chi(1) p^{t_{\chi}} = d(Z(G))p^{b+1}. \]
From Lemma \ref{lemma:m(chi)}, $|X_G| = d(Z(G))$, and thus, we get $\chi(1) p^{t_{\chi}} = p^{b+1}$, for each $\chi \in X_G$. Therefore, if $\chi \in X_G \cap \nl(G)$, we get $p^{t_{\chi}} = p^{b+1}/p^{b} = p$. This proves the claim.\\
Now, we prove that $G\in \mathcal{T}$. On the contrary, suppose that $\ker(\lambda_{\chi})$ is a normal subgroup of $G$, for some $\chi \in X_G \cap \nl(G)$. Then $|\ker(\chi)| = |\ker(\lambda_{\chi})|$. Since $t_{\chi} = 1$, we get  $|G/\ker(\chi)| = p^{b+1}$. From \cite[Corollary 2.30]{I}, we get $\chi(1)^2 \leq |G/ Z(\chi)| \leq |G/ \ker(\chi)|$. Then $p^{2b} \leq p^{b+1}$, which implies that $b \leq 1$. Now if $b =1$, we get $|G/\ker(\chi)| = p^2$, and hence, $G' \subset \ker(\chi)$. This is a contradiction. \qed

\begin{example} \label{example:newclassinT}
\textnormal{\begin{enumerate}
		\item [(i)] 	Consider 
		\[ G = \phi_{3}(2111)e = \langle \alpha, \alpha_{1}, \alpha_{2}, \alpha_{3} ~|~ [\alpha_{i}, \alpha] = \alpha_{i+1}, \alpha^{p} = \alpha_{1}^{p^2} = \alpha_{i+1}^{p} = 1 ~ (i=1,2) \rangle, \]
		a group of order $p^5$ ($p\geq 3$) (see \cite[Section 4.5]{RJ}). Here $Z(G) = \langle \alpha_{1}^{p}, \alpha_{3} \rangle \cong C_p \times C_p$ and $\cd(G) = \{ 1, p \}$. Since $G$ is not a direct product of an abelian and a non-abelian subgroup, from Lemma \ref{thm:mu(G)=c(G)} and \cite[Lemma 14]{SAcyclic}, we get that $2p^2 \leq c(G)$. Further, since $\Core_{G}(\langle \alpha_{1}, \alpha_{2} \rangle) \cap \Core_{G}(\langle \alpha, \alpha_{2}, \alpha_{3} \rangle) = 1$, we get $c(G) = \mu(G) \leq 2p^2$. Therefore, we get $c(G) = 2p^2$. Then from Proposition \ref{prop:newclassinT}, we get $G \in \mathcal{T}$.
		\item [(ii)] Let $G$ be a $p$-group of order $p^6 ~(p \geq 7)$ belonging to the isoclinic family $\Phi_{15}$ (see \cite{Arxiv}). From \cite[Corollary 25]{SA}, $G$ is a VZ $p$-group, $d(Z(G)) = 2$, $\cd(G) = \{ 1, p^2 \}$ and $c(G) = 2p^3$. Further, since a VZ $p$-group is a GVZ $p$-group too, from Theorem \ref{thm:GVZcod}, we get $G\in \mathcal{T}$. 
	\end{enumerate} }
\end{example}

		  In Subsections \ref{subsec:corep}, \ref{subsec:chainofkernels} and \ref{subsec:GVZgroups}, we study various classes of non-abelian $p$-groups belonging in $\mathcal{T}$.
		\subsection{Core-$p$ $p$-groups} \label{subsec:corep}
	A $p$-group $G$ is called core-$p$ $p$-group if $|H/ \Core_{G}(H)| \leq p$, for every subgroup $H$ of $G$. Lennox {\it et al.} \cite{corep1995} started the study of core-$p$ $p$-groups. They proved that if $G$ is a finite core-$p$ $p$-group ($p\geq 3$), then the nilpotency class of $G$ is at most 3. In \cite{corep1997}, Cutolo {\it et al.} proved that a  core-$p$ $p$-group $(p\geq 3)$ has an abelian normal subgroup of index at most $p^2$. 
	In Lemma \ref{lemma:corep}, we mention some preliminary results on core-$p$ $p$-groups that we use in this article.
	
	\begin{lemma} \label{lemma:corep}
		Let $G$ be a non-abelian core-$p$ $p$-group $(p\geq 3)$. Then
		\begin{enumerate}
			\item [(i)] \textnormal{\cite[Corollary 2.2]{corep1995}} $G$ is metabelian.
			\item [(ii)] \textnormal{\cite[Lemma 2.9]{Wilkens}} Let $A$ be an abelian subgroup of $G$ of maximum order. Then $|G:A| \leq p^2$, and if $|G:A| = p^2$, then $|A:Z(G)| = p$.
		\end{enumerate} 
	\end{lemma}

\begin{theorem} \label{thm:kerlambdanotnormal}
	Let $G$ be a non-abelian core-$p$ $p$-group of order $p^n ~(p\geq 3)$. Then
	\begin{enumerate}
		\item [(i)] $\cd(G) = \{ 1, p\}$.
		\item [\rmfamily(ii)]   If $d(Z(G)) = 2$ and $n> 2e+3$, where $\exp(G) = p^e ~(e\geq 1)$, then $G\in \mathcal{T}$.
		\item [\rmfamily(iii)] If $G \in \mathcal{T}$ and $X_G \in \mathcal{S}$, then $c(G) = \sum_{\chi \in X_{G}}\cod(\chi)$. Further, if $Z(G)$ is cyclic, then $c(G) = p^{n-1}$.
	\end{enumerate}
\end{theorem}
	\noindent \emph{Proof.} To prove (i), let $A$ be an abelian subgroup of maximum order in $G$. From Lemma \ref{lemma:corep}(ii), $|G:A| \leq p^2$. If $|G:A| = p$, then $\cd(G) = \{ 1, p \}$ from \cite[Theorem 12.11]{I}. Otherwise, if $|G:A| = p^2$, then $|G:Z(G)| = |G:A||A:Z(G)| = p^3$. Hence $\cd(G) = \{ 1, p \}$.\\
\noindent	  To prove (ii), suppose $d(Z(G)) = 2$ and $n>2e+3$.
	 Let $X_G\in \mathcal{S}$.
	 Since $d(Z(G)) = 2$, we have $|X_{G}| = 2$ (from Lemma \ref{lemma:m(chi)}). Observe that $|X_{G} \cap \lin(G)| = 1$, or $0$. \\
	 {\bf Case I} $(|X_{G} \cap \lin(G)| = 1)$: Suppose $\eta \in X_{G}\cap \lin(G)$. Then $|G/ \ker(\eta)| \leq p^e \Rightarrow |\ker(\eta)| \geq p^{n-e}$. Now, suppose $\chi \in X_{G} \cap \nl(G)$. Then $\chi(1) = p$. Let $H_{\chi} \leq G$ and $\lambda_{\chi} \in \lin(H_{\chi})$ be as defined in Remark \ref{remark:ford}. Since $H_{\chi}/\ker(\lambda_{\chi})$ is cyclic, we get $|H_{\chi}/\ker(\lambda_{\chi})| \leq p^e \Rightarrow |\ker(\lambda_{\chi})| \geq p^{n-e-1}$. Now if possible, suppose that $\ker(\lambda_{\chi})$ is normal in $G$. Then $\ker(\chi) = \Core_{G}(\ker(\lambda_{\chi})) = \ker(\lambda_{\chi})$, which implies that $|\ker(\chi)| \geq p^{n-e-1}$. Since $\chi, \eta \in X_{G}$, we get $|\ker(\chi)\cdot \ker(\eta)| = |\ker(\chi)|\cdot |\ker(\eta)| \geq p^{2n-2e-1}$.
	 Since $n > 2e+3$, we have $2n-2e-1> 2n -2e -3=2n - (2e+3) > 2n - n = n$. Therefore, we get $|\ker(\chi)\cdot \ker(\eta)| > p^n = |G|$, a contradiction.\\
{\bf Case II} $(|X_{G} \cap \lin(G)| = 0)$:	 Let $X_{G} = \{ \chi, \psi \}$, where $\chi(1), \psi(1)>p$. Let $K_{\chi},K_{\psi} \leq G$ and $\lambda_{\chi} \in \lin(K_{\chi})$, $\rho_{\psi} \in \lin(K_{\psi})$ be as defined in Remark \ref{remark:ford}. By the discussion in Case I, we get $|\ker(\lambda_{\chi})|, |\ker(\rho_{\psi})| \geq p^{n-e-1}$. \\
{\bf Claim:} $|\ker(\lambda_{\chi})|$ and $|\ker(\rho_{\psi})|$ are not normal in $G$. 
On the contrary, suppose that one of them, say $\ker(\lambda_{\chi})$, is normal in $G$. Then $|\ker(\chi)| = |\ker(\lambda_{\chi})| \geq p^{n-e-1}$. Since $G$ is a core-$p$ $p$-group, we get $|\ker(\rho_{\psi})| \leq p|\Core_{G}(\ker(\rho_{\psi}))| \Rightarrow |\ker(\psi)| \geq p^{n-e-2}$. Since $\ker(\chi) \cap \ker(\psi) = 1$, we get $|\ker(\chi)\cdot \ker(\psi)| = |\ker(\chi)|\cdot |\ker(\psi)| \geq p^{2n-2e-3}$. Since $n> 2e+3$, we get $|\ker(\chi)\cdot \ker(\psi)| > p^{n} = |G|$, a contradiction. This proves the claim, and hence, $G\in \mathcal{T}$.\\
\noindent	Now we prove (iii). Suppose $X_G\in \mathcal{S}$ and let $\chi \in X_G \cap \nl(G)$. 
	 Let $H_{\chi} \leq G$ and $\lambda_{\chi} \in \lin(H_{\chi})$ be as defined in Remark \ref{remark:ford}.
	Since $G\in \mathcal{T}$, $\ker(\lambda_{\chi})$ is not normal in $G$. Further, since $G$ is a core-$p$ $p$-group, $|\ker(\lambda_{\chi})| = p|\ker(\chi)|$. Then
	\begin{equation*} \label{eq:codnlinear3}
		\left| \frac{H_{\chi}}{\ker(\lambda_{\chi})} \right| = \frac{|H_\chi|}{p|\ker(\chi)|}  = \frac{1}{p\chi(1)}\left| \frac{G}{\ker(\chi)} \right| = \frac{1}{p}\cod(\chi) \neq 1.
	\end{equation*}
	Since $\chi(1) = p$, we get
	\[ 	\frac{p}{p-1}d(\chi) =  \frac{p}{p-1}\chi(1)|\Gamma(\chi)| = \frac{p}{p-1}\chi(1) \phi\left( \left| \frac{H_{\chi}}{\ker(\lambda_{\chi})} \right| \right)= \frac{p}{p-1}\phi\left( \cod(\chi) \right)=  \cod(\chi). \]
	Then from Remark \ref{remark:m(chi)} and Lemma \ref{lemma:codlinearchar}, we get $c(G) = \frac{p}{p-1} \sum_{\chi \in X_{G}}d(\chi) = \sum_{\chi \in X_{G}} \cod(\chi).$ Now, if $Z(G)$ is cyclic, $|X_G| = 1$ (from Lemma \ref{lemma:m(chi)}). Let $X_G = \{ \chi \}$, for some faithful irreducible character $\chi$ of $G$. Since $\chi(1) = p$, we get $c(G) = \cod(\chi) = p^{n-1}$. \qed\\

	In the next subsection, we study the relation between $c(G)$ and codegrees of $\chi \in X_G$ for various classes of maximal class $p$-groups.
	
	\subsection{Groups where kernels of nonlinear irreducible characters constitute a chain} \label{subsec:chainofkernels}
	
	\noindent	In \cite[p. 440]{YBbook}, Berkovich suggested to study those $p$-groups where kernels of irreducible characters form a chain. In \cite[Corollary 5.2]{Lewischain}, Lewis proved that such groups are cyclic, and modified the Berkovich's question to the study of $p$-groups where kernels of nonlinear irreducible characters form a chain.
	In \cite{QianWang2008}, Qian and Wang proved the following result.

	\begin{lemma} \textnormal{\cite[Theorem 2.5]{QianWang2008}} \label{lemma:chainofkernels}
		Let $G$ be a non-abelian $p$-group. Then the following statements are equivalent.
		\begin{enumerate}
			\item $\Kern(G) = \{ \ker(\chi) ~ |~ \chi \in \nl(G) \}$ is a chain with respect to inclusion.
			\item $G$ is one of the following groups:\\
			\textnormal{(2.1)} $G'$ is a unique minimal normal subgroup of $G$.\\
			\textnormal{(2.2)} $G$ is of maximal class.
		\end{enumerate}
	\end{lemma}
	\noindent By \cite[Remark 1.1]{QianWang2008}, if $G$ is a finite non-abelian $p$-group such that $G'$ is a unique minimal normal subgroup of $G$, then $Z(G) \supseteq G'$ is cyclic, $G/Z(G)$ is an elementary abelian group of order $p^{2m}$, and all its nonlinear irreducible characters are faithful and of degree $p^m$. It is easy to see that $G$ is, in fact, a VZ $p$-group with cyclic center. Further, a $p$-group of maximal class has a center of order $p$ (see \cite[Exercise 1, p. 114]{YBbook}). Hence from Lemma \ref{lemma:cycliccenterinT}, we deduce Corollary \ref{cor:chainKernelsinT}.
	\begin{corollary}\label{cor:chainKernelsinT}
		Let $G$ be a non-abelian $p$-group such that $\Kern(G)$ is a chain with respect to inclusion. Then $G\in\mathcal{T}$.
	\end{corollary}
	
 We get $c(G)$ for the groups of type 2.1 of Lemma \ref{lemma:chainofkernels} by using \cite[Corollary 4]{SA}.
	\begin{corollary}
		Suppose $G$ is a finite non-abelian $p$-group such that $G'$ is a unique minimal normal subgroup of $G$. Then $c(G) = |G/Z(G)|^{1/2} |Z(G)|$.
	\end{corollary}
	
	 Now, we deal with the groups of type 2.2 of Lemma \ref{lemma:chainofkernels}, i.e., maximal class $p$-groups. 
	In Lemma \ref{lemma:basicresultsmaximal} and Remark \ref{remark:N<frattini}, we mention few background results on maximal class $p$-groups which we use further.
	\begin{lemma}  \label{lemma:basicresultsmaximal}
		Let $G$ be a $p$-group of maximal class and order $p^n$. Then
		\begin{enumerate}
			\item \textnormal{\cite[Exercise 1, p. 114]{YBbook}} $|G/G'| = p^2$ and $|Z(G)| = p$. If $1\leq i \leq n-2$, then $G$ has only one normal subgroup of order $p^i$.
			\item \textnormal{\cite[Exercise 1, p. 114]{YBbook}} If $p\geq 3$ and $n>3$, then $G$ has no cyclic normal subgroups of order $p^2$.
			\item \textnormal{\cite[Theorem 9.5]{YBbook}} If $n \leq p+1$, then $\Phi(G)$ and $G/Z(G)$ have exponent $p$, where $\Phi(G)$ denotes the Frattini subgroup of $G$.
		\end{enumerate}
	\end{lemma}

	\begin{remark} \label{remark:N<frattini}
		\textnormal{	Let $G$ be a $p$-group of maximal class and order $p^n$.
			\begin{enumerate}
				\item From Lemma \ref{lemma:basicresultsmaximal}(1) and \cite[p. 26]{YBbook}, we get 
				$G' = \Phi(G)$ is of index $p^2$ in $G$. 
				\item If $n\leq p+1$, then from Lemma \ref{lemma:basicresultsmaximal}(3) it is easy to see that $\exp(G') = p$.
				\item Let $N$ be the unique proper normal subgroup of order $p^m$ in  $G$. Then from \cite[Proposition 8.2(a)]{YB2015}, $N \leq \Phi(G) = G'$.
			\end{enumerate}
		}
	\end{remark}
	
 We begin with the following result for normally monomial $p$-groups with cyclic center.
	
	\begin{proposition} \label{prop:codnMcyclic}
		Let $G$ be a normally monomial $p$-group $(p\geq 3)$ with cyclic center. Suppose $X_G = \{ \chi \} \in \mathcal{S}$. If there does not exist any cyclic normal subgroup of index $b(G)$ in $G$, then $c(G) \leq \frac{b(G)}{p} \cod(\chi)$. Moreover, for each $n\geq 3$, there exists a normally monomial $p$-group $(p\geq 3)$ of order $p^n$ with cyclic center for which the equality occurs.
	\end{proposition}
	
	\noindent \emph{Proof.}  Let $X_G = \{\chi\}\in \mathcal{S}$. From \cite[Proposition 3]{AM}, $\chi(1) = b(G) = p^a$ (say) and $\chi$ can be linearly induced from an abelian normal subgroup, say $A$, of index $b(G)$ in $G$. Let $\chi = \lambda\ind_{A}^{G}$, for some $\lambda \in \lin(A)$. Since $A$ is not cyclic, $\ker(\lambda) \neq 1$. Hence, we have $p|\ker(\chi)| \leq |\ker(\lambda)|$. Then 
	\[ \left| \frac{A}{\ker(\lambda)} \right|  \leq \frac{|A|}{p|\ker(\chi)|} = \frac{1}{p\chi(1)} \left| \frac{G}{\ker(\chi)} \right| = \frac{1}{p} \cod(\chi). \]
	Let $\cod(\chi) = p^b$, for some $b\geq 1$. Then $d(\chi) = \chi(1)|\Gamma(\chi)| \leq p^a |\Gamma(\lambda)| = p^a \phi(|A/\ker(\lambda)|) \leq p^{a}\phi(p^{b-1}) = p^{a+b-2}(p-1)$. From
	Remark \ref{remark:m(chi)}, $c(G) = \frac{p}{p-1}d(\chi) \leq p^{a+b-1} = \frac{b(G)}{p} \cdot \cod(\chi)$.\\
	Now, for each $n\geq 3$, consider
	\[ G = \langle a,b ~|~ b^{-1}a^{-1}ba= a^{p^{n-2}}, a^{p^{n-1}} = b^{p} = 1 \rangle, \] a non-abelian $p$-group of order $p^n$ ($p \geq 3$).
	 Here $Z(G) = \langle a^{p} \rangle \cong C_{p^{n-2}}$ (and hence, $G\in \mathcal{T}$) and $\cd(G) = \{ 1, p \}$. It is easy to see that $c(G) = p^{n-1}$. Now, suppose $X_G = \{ \chi \}\in \mathcal{S}$. Here $\chi(1) = p$ and $\ker(\chi) = 1$. Then 
		\[ c(G) = p^{n-1} = \cod(\chi) =   \frac{b(G)}{p}\cod(\chi). \eqno \qed \]

\begin{remark} \label{remark:nMonomial}
	\textnormal{ We do have an example of a $p$-group satisfying the hypothesis of Proposition \ref{prop:codnMcyclic} such that $c(G) < \frac{b(G)}{p} \cod(\chi)$. Consider
		\[ G =  \langle x, y, z: x^{p^2} = y^p = z^p = 1, xy = yx^{p+1}, xz = zxy, yz = zy \rangle, \]
		which is a $p$-group of order $p^4$ ($p \geq 3$) with $Z(G) = \langle x^p \rangle \cong C_p$ and $\cd(G) = \{ 1, p \}$. Here, $b(G) = p$. Suppose $X_G = \{ \chi \} \in \mathcal{S}$. Since $\exp(G) = p^2$, there does not exist any cyclic normal subgroup of index $b(G)$ in $G$. From Example \ref{example:intro2}, we have $c(G) < \cod(\chi) = \frac{b(G)}{p} \cod(\chi)$.
	}
\end{remark}
	
 For normally monomial maximal class $p$-groups of order $p^n$, we deduce the following result with the help of Lemma \ref{lemma:normallymonomial} and Proposition \ref{prop:codnMcyclic}.

	\begin{proposition} \label{prop:nMaximal}
	Let $G$ be a normally monomial maximal class $p$-group of order $p^n~ (p\geq 3, n>3)$. Let $X_G = \{ \chi \} \in \mathcal{S}$. Then $c(G) \leq \frac{b(G)}{p}\cod(\chi)$, and $b(G)\cdot p \leq c(G) \leq b(G) \cdot \exp(A)$, where $A$ is an abelian subgroup of index $b(G)$ in $G$. Further,
	\begin{enumerate}
		\item [\rmfamily(ii)] if $b(G) = p$ with $n \leq p+1$, then $c(G) = p^2$ or $p^3$. 
		\item [\rmfamily(iii)] If $b(G) > p$, then $b(G)\cdot p \leq c(G) \leq b(G)\cdot \exp(G')$. Moreover, if $n\leq p+1$, then $c(G) = b(G) \cdot p$. 
	\end{enumerate}
\end{proposition}
\noindent \emph{Proof.}
	 Suppose $X_G = \{ \chi\} \in \mathcal{S}$, which implies that $\ker(\chi) = 1$ and $\chi(1) = b(G)$ \cite[Proposition 3]{AM}.\\
	 {\bf Claim:} $G$ has no cyclic normal subgroup of index $b(G)$ in $G$. On the contrary, suppose that $G$ has a cyclic normal subgroup $A$ of index $b(G)$ in $G$.  Now, let $b(G) = p^e$. Then $p^{2e} \leq |G|/|Z(G)| = p^{n-1}$. Since $|G/A| = p^e$, we get $|A| \geq p^{(n+1)/2}$. Therefore, for $n>3$, $|A| > p^2$. Then $A$ has a cyclic subgroup, say $B$, of order $p^2$. Here $B$ is a characteristic subgroup of $A$, and hence, is a normal subgroup of $G$. This is a contradiction
	  from Lemma \ref{lemma:basicresultsmaximal}(2). It proves the claim. \\
	   Then from Proposition \ref{prop:codnMcyclic}, we get $c(G) \leq \frac{b(G)}{p}\cod(\chi)$. Now, from \cite[Proposition 3]{AM}, there exists an abelian normal subgroup $A$ of $G$ such that $|G/A| = p^b$ and $\chi = \lambda\ind_{A}^{G}$, for some $\lambda \in \lin(A)$.
	Then from Lemma \ref{lemma:normallymonomial}, we have
	\[ b(G)\cdot p \leq c(G) \leq b(G) \cdot \exp(A). \]
	Now, suppose $b(G) = p$ with $n \leq {p+1}$. 
	From Lemma \ref{lemma:basicresultsmaximal}(3), we have $\exp(G/Z(G)) = p$. Since $Z(G) \cong C_{p}$, we get $\exp(G) \leq p^2 \Rightarrow \exp(A) \leq p^2$. This implies that $c(G) = p^2$ or $p^3$.\\
	\noindent Now, suppose $b(G) > p$. Then $|G/A| > p$. From Lemma \ref{lemma:basicresultsmaximal}(1) and Remark \ref{remark:N<frattini}(3), $A$ is contained in $G'$. Then from Lemma \ref{lemma:normallymonomial}, we get $b(G)\cdot p \leq c(G) \leq b(G) \cdot \exp(A) \leq b(G) \cdot \exp(G')$.
	 Further, let $n \leq {p+1}$. From Remark \ref{remark:N<frattini}(2), we have $\exp(A) = \exp(G') = p$. 
	Therefore, we get $c(G) = b(G) \cdot p$.  \qed
	
	\begin{example} \label{example1:nMMaximal}
		\textnormal{
Let $G$ be a $p$-group $(p\geq 5)$ of order $p^6$ in the isoclinic family $\Phi_{i}$, where $i\in \{ 35, 36, 38 \}$. Then $G$ is a normally monomial maximal class $p$-group such that $\cd(G) = \{ 1,p \}$ if $G\in \Phi_{35}$, and $\cd(G) = \{ 1, p, p^2 \}$ if $G\in \Phi_{i}$ where $i\in \{ 36, 38 \}$  (see \cite[Section 4.1]{RJ}). Then from Proposition \ref{prop:nMaximal}, we get
\[ c(G) = 
\begin{cases}
	p^2 \text{ or } p^3, &\text{ if } G\in \Phi_{35}\\
	p^3, &\text{ if } G\in \Phi_{i} \text{ where } i\in \{ 36, 38 \}.
\end{cases}
 \]
	}
	\end{example}

%
	

	 Suppose $G$ is a maximal class $p$-group of order $p^n$ with $\cd(G) = \{ 1, p^{a} \}$, for some integer $a\geq 1$. From \cite[Theorem 22.12]{YBbook}, $\Phi(G)$ is abelian. Then from Remark \ref{remark:N<frattini}(1), $G'$ is abelian, which makes $G$ a metabelian $p$-group. Since a metabelian $p$-group is a normally monomial $p$-group too, a range of $c(G)$ is given by Proposition \ref{prop:nMaximal}. In Proposition \ref{prop:maximalpgroup1}, we deal with the case when $\cd(G) = \{ 1, p, p^{a} \}$, for some integer $a>1$.
	
	\begin{proposition} \label{prop:maximalpgroup1}
		Suppose $G$ is a maximal class $p$-group $(p\geq 3)$ of order $p^n$ such that $n\leq p+1$ and $\cd(G) = \{ 1, p, p^{b} \}$, for some integer $b>1$. Then $c(G) = p^{b+1}$ or $p^{b+2}$.
	\end{proposition}
	\noindent \emph{Proof.}  Since $Z(G)$ is cyclic, let $X_G = \{ \chi \}\in \mathcal{S}$. From \cite[Lemma 21]{SAcyclic}, $\chi(1) = p^b$. Let $H_{\chi}\leq G$ and $\lambda_{\chi} \in \lin(H_{\chi})$ be as defined in Remark \ref{remark:ford}. From Lemma \ref{lemma:basicresultsmaximal}(3), we get $\exp(G) \leq p^2$. Hence, we get $\frac{H_{\chi}}{\ker(\lambda_{\chi})} \cong C_{p}$ or $C_{p^2}$. This implies that $|\Gamma(\lambda_{\chi})| = \phi(p)$ or $\phi(p^2)$ respectively. Then
	\[ 
	d(\chi) = \chi(1)|\Gamma(\chi)| =  
	\begin{cases}
		p^{b}(p-1), & \text{ if } |\Gamma(\lambda_{\chi})| = \phi(p)\\
		p^{b+1}(p-1), & \text{ if } |\Gamma(\lambda_{\chi})| = \phi(p^2).
	\end{cases}
	 \]
	  Therefore, from Remark \ref{remark:m(chi)}, we get $c(G)\in \{ p^{b+1}, p^{b+2} \}$. \qed

\subsection{GVZ $p$-groups} \label{subsec:GVZgroups}
In this subsection, we study the relation between $c(G)$ and the character codegrees of GVZ $p$-groups and prove that GVZ $p$-groups $(p\geq 3)$ belong to $\mathcal{T}$.
GVZ-groups were first introduced by Nenciu in \cite{Nenciu2}. A group $G$ is called a GVZ-group if $\chi(g) = 0$, for all $g\in G\setminus Z(\chi)$ and all $\chi\in \nl(G)$. Readers can see \cite{BL,Lewischain, Nenciu2, Nenciu} for background results on GVZ-groups. In \cite[Proposition 1.2]{Nenciu2}, Nenciu proved that GVZ-groups are nilpotent. Examples produced in \cite{Lewischain,Nenciu} show that there exist GVZ-groups of arbitrarily
high nilpotency class. 
  Now, when $G$ and $H$ are nilpotent groups, we have $c(G\times H) = c(G) + c(H)$ (see \cite{MG}). Further, since a GVZ-group is nilpotent, it can be expressed as a direct product of its Sylow subgroups. Hence we focus on GVZ $p$-groups in this article. 
We also separately study $\CM_{p-1}$ $p$-groups, which form a sub-class of GVZ $p$-groups.
 We have the following result for GVZ $p$-groups.\\

\noindent {\bf Proof of Theorem \ref{thm:GVZcod}.} Let $G$ be a GVZ $p$-group $(p\geq 3)$ and suppose $\chi \in \Irr(G)$. Let $H_{\chi}\leq G$ and $\lambda_{\chi} \in \lin(H_{\chi})$ be as defined in Remark \ref{remark:ford}. Suppose $|\ker(\lambda_{\chi})|/|\ker(\chi)| = a_{\chi}$. Then
\begin{equation*}
	\left| \frac{H_{\chi}}{\ker(\lambda_{\chi})} \right| = \frac{|H_{\chi}|}{a_{\chi}|\ker(\chi)|} = \frac{|G|}{a_{\chi}\chi(1) |\ker(\chi)|} = \frac{\cod(\chi)}{a_{\chi}},
\end{equation*}
which implies that $|\Gamma(\chi)| = |\Gamma(\lambda_{\chi})| = \phi({\cod(\chi)}/{a_{\chi}})$. Now, since $G$ is a GVZ $p$-group, we have $\chi(g) = 0$, for all $g\in G \setminus Z(\chi)$. Then from \cite[Corollary 2.30]{I}, we get $|G/Z(\chi)| = \chi(1)^2$. From \cite[Lemma 2.27]{I}, we have $\chi\restr_{Z(\chi)} = \chi(1) \mu$, for some $\mu \in \lin(Z(\chi))$. Hence $\ker(\chi) = \ker(\mu)$, and we get
\begin{equation} \label{eq:GVZcod}
	\left| \frac{Z(\chi)}{\ker(\mu)} \right| = \frac{|Z(\chi)|}{|\ker(\chi)|} = \frac{|G|}{\chi(1)^{2} |\ker(\chi)|} = \frac{\cod(\chi)}{\chi(1)}.
\end{equation}
Hence, we get $|\Gamma(\chi)| = \phi(\cod(\chi)/\chi(1))$. Then from the above discussions, we have $|\Gamma(\chi)| = \phi({\cod(\chi)}/{a_{\chi}})  = \phi(\cod(\chi)/\chi(1))$, and thus, $a_{\chi} = \chi(1)$. 
Now, if $\chi \in \nl(G)$, then we get $|\ker(\lambda_{\chi})|/|\ker(\chi)| = a_{\chi} > 1$. Therefore, for each $\chi\in \nl(G)$, we get $\ker(\lambda_{\chi})$ is not normal in $G$ and thus, $G \in \mathcal{T}$.\\
\noindent Now, let $G$ be a GVZ $p$-group ($p$ is any prime) and suppose $X_G \in \mathcal{S}$.
Let $\chi\in X_G\cap \nl(G)$. Then $\chi\restr_{Z(\chi)} = \chi(1) \mu$, for some $\mu \in \lin(Z(\chi))$, and thus from Equation \eqref{eq:GVZcod}, we get
$d(\chi) = \chi(1) |\Gamma(\chi)| =  \chi(1) \phi(\cod(\chi)/ \chi(1)) = \phi(\cod(\chi))$. Then 
\begin{equation*}\label{eq:codGVZnl}
	\frac{p}{p-1}d(\chi) = \cod(\chi).
\end{equation*}
Therefore, from Remark \ref{remark:m(chi)} and Lemma \ref{lemma:codlinearchar}, we get
\begin{equation*}
	c(G) = \frac{p}{p-1}\sum_{\chi \in X_G}d(\chi)= \sum_{\chi \in X_{G}} \cod(\chi). \eqno \qed
\end{equation*}

\begin{corollary} \label{cor:GVZcyclic}
	Let $G$ be a GVZ $p$-group with cyclic center. Then $c(G) = |G/Z(G)|^{1/2} |Z(G)|$.
\end{corollary}
\noindent \emph{Proof.} From Lemma \ref{lemma:m(chi)}, we get $|X_{G}| = 1$, where $X_G \in \mathcal{S}$. Suppose $X_G = \{ \chi \}$. From Theorem \ref{thm:GVZcod}, we get $c(G) = \cod(\chi)$.  Since $\ker(\chi) = 1$, we have $Z(\chi) = Z(G)$  \cite[Lemma 2.27]{I}. Since $G$ is a GVZ $p$-group, from \cite[Corollary 2.30]{I}, we get $\chi(1) = |G/Z(G)|^{1/2}$. Hence
\[ c(G) = \cod(\chi) =  \frac{|G|}{\chi(1)} = \frac{|G|}{|G/Z(G)|^{1/2}} = |G/Z(G)|^{1/2} |Z(G)|. \eqno \qed \]

\noindent In  \cite[Theorem 4.12]{HB}, Behravesh proved that if $G$ is a $p$-group of nilpotency class 2 with cyclic center, then $c(G) = |G/Z(G)|^{1/2} |Z(G)|$. Since every $p$-group of class 2 is a GVZ $p$-group, Corollary \ref{cor:GVZcyclic} generalizes Behravesh's result. \\

\noindent For all $\chi \in \nl(G)$, we have $Z(G) \subset Z(\chi)$. Hence, every VZ $p$-group is a GVZ $p$-group. Thus, the following corollary follows from Theorem \ref{thm:GVZcod}.
\begin{corollary}
	\label{cor:codVZgroup}
	Let $G$ be a VZ $p$-group.
	Then $c(G) = \sum_{\chi \in X_{G}} \cod(\chi). $
\end{corollary}

Now we consider the case of  $\CM_{p-1}$ $p$-groups. A group $G$ is called a $\CM_{p-1}$-group if every normal subgroup of $G$ appears as the kernel of at most $p-1$ irreducible characters of $G$. In \cite{BL}, Burkett and Lewis studied the relation between GVZ $p$-groups and $\CM_{p-1}$ $p$-groups, and proved the following result.
	
	\begin{lemma} \textnormal{\cite[Theorem 8]{BL}} \label{lemma:CMgroup}
		Let $G$ be a $p$-group. Then the following are equivalent.
		\begin{enumerate}
			\item $G$ is a $\CM_{p-1}$-group.
			\item $G$ is a GVZ-group and $|Z(\chi)/\ker(\chi)| = p$ for all $1_{G} \neq \chi \in \Irr(G)$.
			\item $G$ is
			a GVZ-group and every character in $\Irr(G)$ has values in the $p$-th cyclotomic field.
		\end{enumerate}
	\end{lemma}
	
	
	We prove the following result.\\
	
	\noindent {\bf Proof of Theorem \ref{thm:CMgroup}.}
	Let $G$ be a $\CM_{p-1}$ $p$-group. From Lemma \ref{lemma:CMgroup}, $G$ is a GVZ $p$-group and $|Z(\chi)/ \ker(\chi)| = p$, for all $1_{G} \neq \chi \in \Irr(G)$. Then from \cite[Corollary 2.30]{I}, we get $\chi(1)^2 = |G/ Z(\chi)|$. Then 
	\[ \cod(\chi) = \frac{|G/\ker(\chi)|}{\chi(1)} = \frac{|G/\ker(\chi)|\chi(1)}{\chi(1)^2} = \frac{|G|\chi(1)}{|\ker(\chi)||G/Z(\chi)|}= p \chi(1), \forall ~ 1_{G} \neq \chi \in \Irr(G). \]
	Since $\cod(1_{G}) = 1$, we get $\cod(G) = \{ 1, p\chi(1)~ |~ 1_{G} \neq \chi \in \Irr(G) \}$.\\
	 Now, let $\chi \in X_{G} \in \mathcal{S}$.
	From Lemma \ref{lemma:CMgroup}, we get $d(\chi) = \chi(1)|\Gamma(\chi)| = \chi(1) [\mathbb{Q}(\omega_{p}) : \mathbb{Q}] = \chi(1) \phi(p)$. Then $\frac{p}{p-1}d(\chi) = p\chi(1)$. Therefore, from Remark \ref{remark:m(chi)}, we get 
	\[ c(G) = \frac{p}{p-1}\sum_{\chi \in X_{G}} d(\chi) = \sum_{\chi \in X_{G}} p\chi(1). \]
%
	 Since $\cod(\chi) =  p \chi(1)$, for all $1_{G} \neq \chi \in \Irr(G)$, we get $c(G) = \sum_{\chi \in X_{G}} \cod(\chi)$. \\
\noindent Now, to prove (iii), suppose $\cd(G) = \{ 1 = d_{0}, d_{1}, \ldots, d_{t} \}$. Let $X_G, Y_G \in \mathcal{S}$ such that $Y_{G} \neq X_{G}$. From Lemma \ref{lemma:m(chi)}, we have $|Y_{G}| = |X_{G}| = d(Z(G))$. Since $c(G) =   p \sum_{\nu \in Y_{G}} \nu(1) = p\sum_{\chi \in X_{G}} \chi(1)$, we must have $|\Irr_{d_{i}}(G) \cap Y_{G}| = |\Irr_{d_{i}}(G) \cap X_{G}|$, for each $i$ ($ 0\leq i \leq t$). \qed \\
	
	\noindent One should note that Theorem \ref{thm:CMgroup}(iii) may not hold true for all GVZ $p$-groups. That is, if $G$ is a GVZ $p$-group, then there may exist  $X_{G}, Y_{G} \in \mathcal{S}$ such that $|\Irr_{d_{i}}(G) \cap X_{G}| \neq |\Irr_{d_{i}}(G) \cap Y_{G}|$ for some $i$ $(0\leq i \leq t)$. We consider a group of order $5^6$ in Example \ref{example:differentsets} to illustrate one such instance. 
	\begin{example} \label{example:differentsets}
		\textnormal{Consider \[ G = \langle a, b ~|~ a^{-1}b^{-1}ab = a^5, a^{-5}b^{-1}a^{5}b= a^{5^2}, a^{5^3} = b^{5^3} = 1 \rangle, \] which is a metacyclic group of order $5^6$. Through computation in {\sc Magma}, we get that $G$ is also a GVZ group.
			Here $\cd(G) = \{ 1, 5, 5^2 \}$, $Z(G) = \langle a^{5^2}, b^{5^2} \rangle \cong C_{5} \times C_{5}$, $G' = \langle a^{5} \rangle \cong C_{5^2}$ and $G/G' = \langle aG', bG' \rangle \cong C_{5}\times C_{5^3}$. Then from \cite[Lemma 3.10]{BehraveshMetacyclic}, we get $c(G) = 2\cdot 5^3=250$. Now we find a minimal faithful quasi-permutation representation of $G$ satisfying \eqref{eq:X_G}. Let $\psi \in \lin(G/G')$ given by $\psi = 1_{\langle aG' \rangle} \cdot \psi_{\langle bG' \rangle}$, where $1_{\langle aG' \rangle}$ is the trivial character of $\langle aG' \rangle$ and $\psi_{\langle bG' \rangle}$ is a faithful linear character of $\langle bG' \rangle$. Then $\ker(\psi) = \langle a \rangle$. Now consider $H = \langle a^{5^2}, b \rangle \cong C_{5} \times C_{5^3}$ which is a subgroup of $G$ of order $5^4$. Let $\lambda \in \lin(H)$ given by $\lambda = \lambda_{\langle a^{5^2} \rangle} \cdot 1_{\langle b \rangle}$, where $\lambda_{\langle a^{5^2} \rangle}$ is a faithful linear character of $\langle a^{5^2} \rangle$ and $1_{\langle b \rangle}$ is  the trivial character of $\langle b \rangle$. From Mackey's irreducibility criterion, we get $\chi = \lambda\ind_{H}^{G} \in \Irr(G)$. Here $\ker(\chi) = \Core_{G}(\ker(\lambda)) \subseteq \langle b \rangle$.
			Therefore, we have $\ker(\psi) \cap \ker(\chi) = 1$. Now, suppose $\xi =  \sum_{\sigma \in \Gamma(\psi)}\psi^{\sigma} +  \sum_{\sigma \in \Gamma(\chi)} \chi^{\sigma}$. Then $\xi(1) = \sum_{\sigma \in \Gamma(\psi)}\psi^{\sigma}(1) +  \sum_{\sigma \in \Gamma(\chi)} \chi^{\sigma}(1) = \psi(1)|\Gamma(\psi)| + \chi(1)|\Gamma(\chi)|$. Here $|\Gamma(\psi)| = \phi(|G/\ker(\psi)|)$ and $|\Gamma(\chi)| = |\Gamma(\lambda\ind_{H}^{G})| \leq |\Gamma(\lambda)| = \phi(|H/ \ker(\lambda)|) = \phi(5)$. Since $|\Gamma(\lambda\ind_{H}^{G})| \geq \phi(5)$, we get $|\Gamma(\lambda\ind_{H}^{G})| = \phi(5)$. Thus, $\xi(1) = \phi(5^3) + 5^2 \phi(5) = 2 \cdot 5^2(5-1)$. From \cite[Lemma 4.5]{HB}, we get $m(\psi) = 5^2$ and $m(\lambda) = 1$. Through routine computation, we get $m(\chi) = 5^2m(\lambda) = 5^2$, and hence, $m(\xi) = 2\cdot 5^2$. Then $\xi(1) + m(\xi) = 2\cdot 5^3 = c(G)$, and therefore, $X = \{ \psi, \chi \} \in \mathcal{S}$.}
		
			\textnormal{Now, consider $K = \langle a, b^5 \rangle$. Then $K$ is a normal subgroup of index $5$ in $G$. Here $K' = \langle a^{5^2} \rangle \cong C_{5}$ and $K/K' = \langle aK', b^{5} K' \rangle \cong C_{5^2} \times C_{5^2}$. Let $\eta \in \lin(K/K')$ given by $\eta = 1_{\langle aK' \rangle} \cdot \eta_{\langle b^5 K'\rangle}$, where $1_{\langle aK' \rangle}$ is the trivial character of $\langle aK' \rangle$ and $\eta_{\langle b^{5}K' \rangle}$ is a faithful linear character of $\langle b^{5} K' \rangle$. Let $\rho = \eta\ind_{K}^{G}$. Since $\inertiagroup_{G}(\eta) = K$, we get $\rho \in \Irr(G)$. Here $\ker(\rho) = \Core_{G}(\ker(\eta)) \subseteq \langle a \rangle$. Then we have $\ker(\rho) \cap \ker(\chi) = 1$. Suppose $\xi' = \sum_{\sigma \in \Gamma(\rho)}\rho^{\sigma} +  \sum_{\sigma \in \Gamma(\chi)} \chi^{\sigma}$. Proceeding as above, we get $\xi'(1) = \rho(1)|\Gamma(\rho)| + \chi(1) |\Gamma(\chi)| \leq 5|\Gamma(\eta)| + 5^2|\Gamma(\lambda)| = 5\phi(|K/\ker(\eta)|) + 5^2 \phi(|H/\ker(\lambda)|) = 5\phi(5^2) + 5^2\phi(5) = 5^2(5-1)$. From \cite[Lemma 4.5]{HB}, we get $m(\eta) = 5$. Through routine computation, we get $m(\rho) = 5m(\eta) = 5^2$, and hence, $m(\xi') = 2\cdot 5^2$. Hence $\xi'(1) + m(\xi') \leq 2\cdot 5^3$. Since $c(G) = 2\cdot 5^3$, we must have $\xi'(1) + m(\xi') = 2\cdot 5^3$. Therefore, $Y = \{ \rho, \chi \} \in \mathcal{S}$.}
			
			\textnormal{Now let $d_{0} = 1$, $d_{1} = 5$ and $d_{2} = 5^2$ and suppose $n_{i} = |\Irr_{d_{i}}(G) \cap X|$, for each $0\leq i\leq 2$. Then $n_{0} = 1$, $n_{1} = 0$ and $n_{2}=1$. This shows that $|\Irr_{d_{i}}(G) \cap Y| \neq n_{i}$, for $i=0,1$. }
	\end{example}

	\section{Results and Discussions} \label{sec:results}

		From Theorem \ref{lemma:generalcodegreerelation}, we deduce that if $G$ is a non-abelian $p$-group ($p\geq 3$) and $X_G \in \mathcal{S}$, then
	\begin{equation} \label{eq:introgeneralpgroup}
		c(G) \leq \sum_{\chi \in X_{G}}\chi(1)\cod(\chi).
	\end{equation} 
We construct a class $\mathcal{T}$ of non-abelian $p$-groups and obtain an improved bound of $c(G)$ in terms of character codegrees for groups belonging to $\mathcal{T}$. 
	In the course of this paper, we obtain various classes of non-abelian $p$-groups which belong to $\mathcal{T}$. Through computation in {\sc Magma}, we observe that all the groups of order $p^5$ $(p=5,7)$ belong to $\mathcal{T}$ as well. It is, then, natural to  ask the following question.\\
	
	\noindent {\bf Question 1:} Classify all the non-abelian $p$-groups ($p \geq 3$) which belong to $\mathcal{T}$.\\
	
 For $G\in \mathcal{T}$, we prove that $c(G) = \sum_{\chi \in X_{G}} \cod(\chi)$ if $G$ is a GVZ $p$-group (see Theorem \ref{thm:GVZcod}), or a core-$p$ $p$-group (see Theorem \ref{thm:kerlambdanotnormal}). 	We also have an example of a non-abelian $p$-group for which $c(G) < \sum_{\chi \in X_{G}} \cod(\chi)$ (see Example \ref{example:intro2}). Hence, the following question arises.\\
 
	\noindent {\bf Question 2.} Classify all the non-abelian $p$-groups ($p \geq 3$) in $\mathcal{T}$ such that
	\[ c(G) = \sum_{\chi \in X_{G}} \cod(\chi), \text{ for each } X_G \in \mathcal{S}. \]

	Now consider 
	\[ G = \langle \alpha_{1}, \alpha_{2}, \beta, \beta_1, \beta_2 ~|~ [\alpha_{1}, \alpha_2] = \beta, [\beta, \alpha_{i}] = \beta_i = \alpha_{i}^{5}, \beta^5 = \beta_{i}^5 = 1 ~ (i=1,2) \rangle, \]
	a non-abelian group of order $5^5$ with $d(Z(G)) = 2$ and $\cd(G) = \{ 1, 5 \}$ (see \cite[Section 4.5]{RJ}). Since $Z(G) \subset G'$, from \cite[Theorem 2]{SAcyclic}, we get $X_G \cap \lin(G) = \emptyset$ for all $X_G \in \mathcal{S}$. Through computation in {\sc Magma}, it is easy to see that $G\in \mathcal{T}$, $c(G) = 2\cdot 5^2 = 50$, and for all $X_G \in \mathcal{S}$,  $\cod(\chi) = 5^3$, for each $\chi\in X_G$. Then $c(G)$ can be expressed in the following way.
	\[ c(G) = \frac{1}{5}\sum_{\chi \in X_{G} \cap \nl(G)} \cod(\chi) = \sum_{\chi \in X_{G} \cap \lin(G)}\cod(\chi) + \frac{1}{5}\sum_{\chi \in X_{G} \cap \nl(G)} \cod(\chi), \text{ for each } X_G \in \mathcal{S}. \]
	 On the other hand, consider 
	\[ G = \langle \alpha, \alpha_{1}, \alpha_{2}, \alpha_{3} ~|~ [\alpha_{i}, \alpha] = \alpha_{i+1}, \alpha^{5} = \alpha_{1}^{5^2} = \alpha_{i+1}^{5} = 1 ~ (i=1,2) \rangle, \]
	a non-abelian group of order $5^5$ with $d(Z(G)) = 2$. From Example \ref{example:newclassinT}(i), $G\in \mathcal{T}$ and $c(G) = 2\cdot 5^2 = 50$. Through computation in {\sc Magma}, it is easy to check that there exist $X_G, Y_G \in \mathcal{S}$
	such that $X_G = \{\psi_1, \chi_1\}$, $Y_G = \{ \chi_2, \chi_3 \}$ with $\psi_1\in \lin(G)$, $\chi_{i}\in \nl(G)$ for $1\leq i \leq 3$, $\cod(\psi_1) = \cod(\chi_{2}) = 5^2$ and $\cod(\chi_{j}) = 5^3$ ($j = 1,3$). Then we can express $c(G)$ in the following way.
	\[ c(G) = \cod(\psi_1) + \frac{1}{5}\cod(\chi_1), \text{ which is also equal to } \cod(\chi_2) + \frac{1}{5}\cod(\chi_3). \]
	Hence, in view of Theorem \ref{lemma:generalcodegreerelation}, we have the following question.\\
	
	\noindent {\bf Question 3.} Classify all the non-abelian $p$-groups ($p \geq 3$) in $\mathcal{T}$ such that
	 \[ c(G) = \sum_{\chi \in X_{G}\cap \lin(G)} \cod(\chi) +  p^{\beta} \left( \sum_{\chi \in X_{G}\cap \nl(G)}\cod(\chi) \right) , \text{ for each } X_G \in \mathcal{S} \text{ and for some fixed } \beta \in \mathbb{Z}.  \]


%
%
%
%

	\section{Acknowledgements}
	Ayush Udeep acknowledges University Grants Commission, Government of India (File: Nov2017-434175) for financial support. The corresponding author acknowledges SERB, Government of India for financial support through grant (MTR/2019/000118). The authors acknowledge University of Auckland for providing remote access to their computational facilities.

\end{document}